\documentclass[12pt]{article}
\usepackage{amssymb}
\usepackage{amsmath}
\newcommand{\lap}{\mbox{$\bigtriangleup$}}

\newcommand{\ra}{{\mbox{$\rightarrow$}}}
\newcommand{\be}{\begin{equation}}
\newcommand{\ee}{\end{equation}}

\newtheorem{mthm}{Theorem}

\newtheorem{thm}{Theorem}[section]

\newtheorem{lem}{Lemma}[section]

\begin{document}

\title{A Liouville theorem for $\alpha$-harmonic functions in $\mathbb{R}^n_+$}

\author{Wenxiong Chen  \quad  Congming Li \quad Lizhi Zhang \quad Tingzhi Cheng  }

\date{\today}
\maketitle
\begin{abstract} In this paper, we consider $\alpha$-harmonic functions in the half space $\mathbb{R}^n_+$:
\begin{equation}
\left\{\begin{array}{ll}
(-\lap)^{\alpha/2} u(x)=0,~u(x)>0, & \qquad x\in\mathbb{R}^n_+, \\
u(x)\equiv0, & \qquad x\notin\mathbb{R}^{n}_{+}.
\end{array}\right.
\label{1}
\end{equation}
We prove that all solutions of (\ref{1}) have to assume the form
 \begin{equation}
u(x)=\left\{\begin{array}{ll}Cx_n^{\alpha/2}, & \qquad x\in\mathbb{R}^n_+, \\
0, & \qquad x\notin\mathbb{R}^{n}_{+},
\end{array}\right.
\label{2}
\end{equation}
for some positive constant $C$.
\end{abstract}
\bigskip

\bigskip
\\
\\
\\

\section{Introduction}

The fractional Laplacian in $R^n$ is a nonlocal pseudo-differential operator, assuming the form
\begin{equation}
(-\Delta)^{\alpha/2} u(x) = C_{n,\alpha} \, \lim_{\epsilon \ra 0} \int_{\mathbb{R}^n\setminus B_{\epsilon}(x)} \frac{u(x)-u(z)}{|x-z|^{n+\alpha}} dz,
\label{Ad7}
\end{equation}
 where $\alpha$ is any real number between $0$ and $2$.
This operator is well defined in $\cal{S}$, the Schwartz space of rapidly decreasing $C^{\infty}$
functions in $\mathbb{R}^n$. In this space, it can also be equivalently defined in terms of the Fourier transform
$$ \widehat{(-\Delta)^{\alpha/2} u} (\xi) = |\xi|^{\alpha} \hat{u} (\xi), $$
where $\hat{u}$ is the Fourier transform of $u$. One can extend this operator to a wider space of distributions.

Let
$$L_{\alpha}=\{u: \mathbb{R}^n\rightarrow \mathbb{R} \mid \int_{\mathbb{R}^n}\frac{|u(x)|}{1+|x|^{n+\alpha}} \, d x <\infty\}.$$

Then in this space, we defined $(-\Delta)^{\alpha/2} u$ as a distribution by
$$< (-\Delta)^{\alpha/2}u(x), \phi> \, = \, \int_{\mathbb{R}^n}  u(x) (-\Delta)^{\alpha/2} \phi (x) dx ,  \;\;\; \forall \, \phi \in C_0^{\infty}(\mathbb{R}^n) . $$

Let
$$\mathbb{R}^n_+ = \{x=(x_1, \cdots, x_n\ \mid x_n > 0 \}$$
be the upper half space. We say that $u$ is $\alpha$-harmonic in the upper half space if
$$\int_{\mathbb{R}^n}  u(x) (-\Delta)^{\alpha/2} \phi (x) dx = 0,  \;\;\; \forall \, \phi \in C_0^{\infty}(\mathbb{R}^n_+) . $$

In this paper, we consider the Dirichlet problem for $\alpha$-harmonic functions
\begin{equation}
\left\{\begin{array}{ll}
(-\lap)^{\alpha/2} u(x)=0,~u(x)>0, & \qquad x\in\mathbb{R}^n_+, \\
u(x)\equiv0, & \qquad x\notin\mathbb{R}^{n}_{+},
\end{array}\right.
\label{1.1}
\end{equation}
It is well-known that
$$
u(x)=\left\{\begin{array}{ll}Cx_n^{\alpha/2}, & \qquad x\in\mathbb{R}^n_+, \\
0, & \qquad x\notin\mathbb{R}^{n}_{+},
\end{array}\right.
$$
is a family of solutions for problem (\ref{1.1}) with any positive constant $C$.

A natural question is: {\em Are there any other solutions?}

Our main objective here is to answer this question and prove
\begin{mthm} Let $0<\alpha<2$, $u\in L_{\alpha}$. Assume $u$
is a solution of
\begin{equation}
\left\{\begin{array}{ll}
(-\lap)^{\alpha/2} u(x)=0,~u(x)>0, & \qquad x\in\mathbb{R}^n_+, \\
u(x)\equiv0, & \qquad x\notin\mathbb{R}^{n}_{+}.
\end{array}\right.
\end{equation}
then
\begin{equation}
u(x)=\left\{\begin{array}{ll}Cx_n^{\alpha/2}, & \qquad x\in\mathbb{R}^n_+, \\
0, & \qquad x\notin\mathbb{R}^{n}_{+},
\end{array}\right.
\end{equation}
for some positive constant $C$.
\label{mthm1}
\end{mthm}

We will prove this theorem in the next section.

\section{The Proof of the Liouville Theorem}

In this section, we prove Theorem \ref{mthm1}. The main ideas are the following.

We first obtain the Poisson representation of the solutions. We show that for $|x-x_r|<r$

\begin{equation}
u(x)=\int_{|y-x_r|>r}P_r(x-x_r,y-x_r)u(y)dy,
\label{P}
\end{equation}
where ~$x_r=(0, \cdots,0,r)$~,~and~$P_r(x-x_r,y-x_r)$~is the Poisson kernel for~$|x-x_r|<r$~:
\begin{eqnarray*}
&&P_{r}(x-x_r,y-x_r)\nonumber\\&=&\left\{\begin{array}{ll}
\frac{\Gamma(n/2)}{\pi^{\frac{n}{2}+1}} \sin\frac{\pi\alpha}{2}
\left[\frac{r^{2}-|x-x_r|^{2}}{|y-x_r|^{2}-r^{2}}\right]^{\frac{\alpha}{2}}\frac{1}{|x-y|^{n}},\qquad& |y-x_r|>r,\\
0,& \text{elsewhere}.
\end{array}
\right.
\end{eqnarray*}

Then, for each fixed $x\in\mathbb{R}^n_+$,  we evaluate first derivatives of $u$ by using (\ref{P}).  Letting $r \rightarrow \infty$,  we derive
$$
\frac{\partial u}{\partial x_i}(x)=0,~~~~~i=1,2,\cdots,n-1.$$
and
$$
\frac{\partial u}{\partial x_n}(x)=\frac{\alpha}{2x_n}u(x).
$$
These yield the desired results.

In the following, we use $C$ to denote various positive constants.
\medskip

{\em Step 1.}

In this step, we obtain the Poisson representation (\ref{P}) for the solutions of (\ref{1.1}).

Let
\begin{equation}
\hat{u}(x)=\left\{
\begin{array}{ll}
\int_{|y-x_r|>r}P_r(x-x_r,y-x_r)u(y)dy, & \qquad |x-x_r|<r, \\
u(x), & \qquad |x-x_r|\geq r.
\end{array}
\right.
\label{max}
\end{equation}
We will prove that ~$\hat{u}$~is ~$\alpha$~-harmonic in ~$B_r(x_r) $. The proof is similar to that in \cite{CL}. It is quite long and complex, hence for reader's convenience, we will present it in the next section.

Let~$w(x)=u-\hat{u}$~,~then
\begin{equation}
\left\{
\begin{array}{ll}
(-\lap)^{\alpha/2}w(x)=0, & \qquad |x-x_r|<r, \\
w(x)\equiv0, & \qquad |x-x_r|\geq r.
\end{array}
\right.
\end{equation}

To show that ~$w\equiv0$~,we employ the following Maximum Principle.

\begin{lem}
(Silvestre, \cite{Si})
 \quad Let $\Omega$ be a bounded domain in $\mathbb{R}^{n}$, and assume that $v$ is a lower semi-continuous function on $\overline{\Omega}$ satisfying
\begin{equation}
\left\{\begin{array}{ll}
(-\lap)^{\frac{\alpha}{2}}v\geq0,  &\mbox{in } \Omega,\\
v\geq0,&\mbox{on } \mathbb{R}^{n}\backslash\Omega.
\end{array}
\right.
\end{equation}
then $v\geq0$ in $\Omega$.
\end{lem}

Applying this lamma to both ~$v=w$~ and ~$v=-w$~,~we conclude that
$$
w(x)\equiv0.
$$
Hence
$$
\hat{u}(x)\equiv u(x).
$$
This verifies (\ref{P}).
\medskip

{\em Step 2.}

We will show that for each fixed ~$x\in\mathbb{R}^n_+$~,
\begin{equation}
\frac{\partial u}{\partial x_i}(x)=0,~~~~~i=1,2,\cdots,n-1.
\label{b1}
\end{equation}
and
\begin{equation}
\frac{\partial u}{\partial x_n}(x)=\frac{\alpha}{2x_n}u(x).
\label{b2}
\end{equation}

From (\ref{b1}), we conclude that ~$u(x)=u(x_n)$, and this, together with (\ref{b2}), immediately implies
\begin{equation}
u(x)=Cx_n^{\alpha/2},
\end{equation}
therefore
\begin{equation}
u(x)=\left\{
\begin{array}{ll}
Cx_n^{\alpha/2}, & \qquad x\in\mathbb{R}^n_+, \\
0, & \qquad x\notin\mathbb{R}^n_+.
\end{array}
\right.
\end{equation}
And this is what we want to derive.

Now, what left is to prove (\ref{b1}) and (\ref{b2}). Through an elementary calculation, one can derive that, for $i=1,2,\cdots,n-1,$
\begin{eqnarray}
\frac{\partial u}{\partial x_i}(x)&=&\int_{|y-x_r|>r}\left(\frac{-\alpha x_i}{r^2-|x-x_r|^2}+\frac{n(y_i-x_i)}{|y-x|^2}\right)P_r(x-x_r,y-x_r)u(y)dy\nonumber\\
&=&\int_{|y-x_r|>r}\frac{-\alpha x_i}{r^2-|x-x_r|^2}P_r(x-x_r,y-x_r)u(y)dy\nonumber\\&\quad&+\int_{|y-x_r|>r}\frac{n(y_i-x_i)}{|y-x|^2}P_r(x-x_r,y-x_r)u(y)dy\nonumber\\
&:=&I_1+I_2.\label{b3}
\end{eqnarray}

For each fixed ~$x\in B_r(x_r)\subset\mathbb{R}^n_+$ and for any given ~$\epsilon >0$~,~we have
\begin{eqnarray}
|I_1|&=&\left|\int_{|y-x_r|>r}\frac{-\alpha x_i}{r^2-|x-x_r|^2}P_r(x-x_r,y-x_r)u(y)dy\right|\nonumber\\
&\leq&\int_{|y-x_r|>r}\left|\frac{-\alpha x_i}{r^2-|x-x_r|^2}\right|P_r(x-x_r,y-x_r)u(y)dy\nonumber\\
&=&\left|\frac{\alpha x_i}{r^2-|x-x_r|^2}\right|\int_{|y-x_r|>r}P_r(x-x_r,y-x_r)u(y)dy\nonumber\\
&=&\left|\frac{\alpha x_i}{r^2-|x-x_r|^2}\right|u(x)\nonumber\\
&\leq&\frac{C}{r}\nonumber\\
&<&\epsilon,~~~\text{for sufficiently large ~$r$~.}\label{b0}
\end{eqnarray}

 Here and in below, the letter ~$C$~ stands for various positive constants.

It is more delicate to estimate $I_2$. For each $R>0$, we divide the region $|y-x_r|>r$ into two parts:
 one inside the ball $|y|<R$ and one outside the ball.

\begin{eqnarray}
|I_2|&=&\left|\int_{|y-x_r|>r}\frac{n(y_i-x_i)}{|y-x|^2}P_r(x-x_r,y-x_r)u(y)dy\right|\nonumber\\
&\leq&\int_{|y-x_r|>r}\left|\frac{n(y_i-x_i)}{|y-x|^2}\right|P_r(x-x_r,y-x_r)u(y)dy\nonumber\\
&=&\int_{\substack{|y-x_r|>r\\|y|>R}}\left|\frac{n(y_i-x_i)}{|y-x|^2}\right|P_r(x-x_r,y-x_r)u(y)dy\nonumber\\&\quad&+\int_{\substack{|y-x_r|>r\\|y|\leq R}}\left|\frac{n(y_i-x_i)}{|y-x|^2}\right|P_r(x-x_r,y-x_r)u(y)dy\nonumber\\
&:=&I_{21}+I_{22}.\label{b4}
\end{eqnarray}

For the ~$\epsilon>0$~ given above, when ~$|y|>R$~, we can easily derive
$$
\left|\frac{n(y_i-x_i)}{|y-x|^2}\right|\leq\frac{n}{|y-x|}\leq\frac{C}{R}<\epsilon,
$$
for sufficiently large $R$. Fix this $R$, then
\begin{eqnarray}
I_{21}&<&\int_{\substack{|y-x_r|>r\\|y|>R}}\epsilon P_r(x-x_r,y-x_r)u(y)dy\nonumber\\
&\leq&\epsilon\int_{|y-x_r|>r}P_r(x-x_r,y-x_r)u(y)dy\nonumber\\
&=&\epsilon u(x)\nonumber\\
&\leq&C\epsilon.\label{b5}
\end{eqnarray}

To estimate $I_{22}$, we employ the expression of the Poisson kernel.
\begin{eqnarray}
I_{22}&=& C \int_{\substack{|y-x_r|>r\\|y|\leq R}}\left[\frac{r^2-|x-x_r|^2}{|y-x_r|^2-r^2}\right]^{\alpha/2}\frac{u(y)}{|x-y|^n}\left|\frac{n(y_i-x_i)}{|y-x|^2}\right|dy\nonumber\\
&\leq& C \int_{\substack{|y-x_r|>r\\|y|\leq R}}\left[\frac{2x_nr-|x|^2}{2y_nr-|y|^2}\right]^{\alpha/2}\frac{u(y)}{|x-y|^n}\frac{1}{|y-x|}dy\nonumber\\
&\leq&C\int_{\substack{|y-x_r|>r\\|y|\leq R}}\left[\frac{2x_nr-|x|^2}{2y_nr-|y|^2}\right]^{\alpha/2}\frac{u(y)}{|x-y|^{n+1}}dy\nonumber\\
&=&C\int_{\substack{|y-x_r|>r\\|y|\leq R, \, y_n>0}}\left[\frac{2x_nr-|x|^2}{2y_nr-|y|^2}\right]^{\alpha/2}\frac{u(y)}{|x-y|^{n+1}}dy\nonumber\\
&\leq&C_R \int_{\substack{|y-x_r|>r\\|y|\leq R, \, y_n>0}}\left[\frac{2x_nr-|x|^2}{2y_nr-|y|^2}\right]^{\alpha/2}\frac{1}{|x-y|^{n+1}}dy.\label{b6}
\end{eqnarray}

Here we have used the fact that the $\alpha$-harmonic function $u$ is bounded in the region
$$D_{R,r}=\left\{y=(y',y_n)\left|\right.|y-x_r|>r,~|y|<R,~y_n>0\right\}.$$
The bound depends on $R$, however is independent of $r$, since $D_{R, r_1} \subset D_{R, r_2}$ for $r_1>r_2$. For each such fixed open domain $D_{R,r}$, the bound of the $\alpha$-harmonic function $u$ can be derived from the interior smoothness ( see, for instance \cite{BKN} and \cite{FW}) and the estimate up to the boundary ( see \cite{RS}).

 Set ~$y=(y',y_n), \sigma=|y'|$, for fixed ~$x$~ and sufficiently large ~$r$, we have

\begin{eqnarray}
\left[\frac{2x_nr-|x|^2}{2y_nr-|y|^2}\right]^{\alpha/2}&\leq&\frac{Cr^{\alpha/2}}{|2y_nr-|y|^2|^{\alpha/2}}\nonumber\\&=&\frac{Cr^{\alpha/2}}{|\sigma^2-2y_nr+y_n^2|^{\alpha/2}}\nonumber\\&=&\frac{Cr^{\alpha/2}}{|(y_n-r)^2+\sigma^2-r^2|^{\alpha/2}}.\label{b7}
\end{eqnarray}
and
\begin{equation}
\frac{1}{|x-y|^{n+1}}\leq\frac{C}{(1+|y|)^{n+1}}\leq\frac{C}{(1+|y'|)^{n+1}}=\frac{C}{(1+\sigma)^{n+1}}.\label{b8}
\end{equation}

For convenience of estimate, we amplify the domain $D_{R,\,r}$ a little bit. 
Define $$\hat{D}_{R,\,r}=\left\{y=(y',y_n)\in\mathbb{R}^n_+\left|\right. |y-x_r|>r,~|y'|\leq R, 0<y_n<\bar{y}_n \right\}.$$  
Here $\bar{y}_n$~ satisfies
\begin{equation}
(\bar{y}_n-r)^2+\sigma^2-r^2=0, \label{b9}
\end{equation}
so that $\bar{y}=(y',\bar{y_n})\in\partial\hat{D}_{R,\,r}\cap\partial B_r(x_r)$. Then it is easy to see that
\begin{equation}
D_{R,\,r}\setminus \hat{D}_{R,\,r}
\label{tian}
\end{equation}

From (\ref{b9}), for sufficiently large ~$r$~(much larger than ~$R$~), we have
$$
\bar{y}_n= r - \sqrt{r^2-\sigma^2}.
$$
Set
\begin{equation}
y_n=r-s\sqrt{r^2-\sigma^2}.\label{b10}
\end{equation}
Then for ~$0<y_n<\bar{y}_n$, 
\begin{equation}
1<s<\frac{r}{\sqrt{r^2-\sigma^2}},\label{b11}
\end{equation}
and
\begin{equation}
dy_n=-\sqrt{r^2-\sigma^2}ds.\label{b12}
\end{equation}

Continuing from the right side of (\ref{b6}), we integrate in the
direction of ~$y_n$~ first, and then integrate with respect to  $y'$, setting $r$ sufficiently large (much larger than $R$), by(\ref{b6}), (\ref{b7}), (\ref{b8}), (\ref{tian}),( \ref{b10}), (\ref{b11}), and(\ref{b12}), we derive
\begin{eqnarray}
I_{22}&\leq&C\int_{\hat{D}_{R,\,r}}\left[\frac{2x_nr-|x|^2}{2y_nr-|y|^2}\right]^{\alpha/2}\frac{1}{|x-y|^{n+1}}dy\nonumber\\
&\leq&C\int_{|y'|<R}\int_0^{\bar{y}_n}\frac{r^{\alpha/2}}{|(y_n-r)^2+\sigma^2-r^2|^{\alpha/2}}dy_n\frac{1}{(1+\sigma)^{n+1}}dy'\nonumber\\
&\leq&C\int_0^R\int_0^{\bar{y}_n}\frac{r^{\alpha/2}}{|(y_n-r)^2+\sigma^2-r^2|^{\alpha/2}}dy_n\frac{1}{(1+\sigma)^{n+1}}\sigma^{n-2}d\sigma \qquad \label{b13}\\
&\leq&C\int_0^R\int_1^{\frac{r}{\sqrt{r^2-\sigma^2}}}\frac{r^{\alpha/2}}{[(s\sqrt{r^2-\sigma^2})^2-(r^2-\sigma^2)]^{\alpha/2}}\left|\sqrt{r^2-\sigma^2}\right|ds\frac{\sigma^{n-2}}{(1+\sigma)^{n+1}}d\sigma\nonumber\\
&=&C\int_0^Rr^{\alpha/2}(r^2-\sigma^2)^{\frac{1-\alpha}{2}}\int_1^{\frac{r}{\sqrt{r^2-\sigma^2}}}\frac{1}{(s^2-1)^{\alpha/2}}ds\frac{\sigma^{n-2}}{(1+\sigma)^{n+1}}d\sigma\nonumber\\
&\leq&C\int_0^Rr^{\alpha/2}(r^2-\sigma^2)^{\frac{1-\alpha}{2}}\int_1^{\frac{r}{\sqrt{r^2-\sigma^2}}}\frac{1}{(s-1)^{\alpha/2}}ds\frac{\sigma^{n-2}}{(1+\sigma)^{n+1}}d\sigma\label{b14}\\
&\leq&C\int_0^Rr^{\alpha/2}r^{1-\alpha}\int_1^{\frac{r}{\sqrt{r^2-\sigma^2}}}\frac{1}{(s-1)^{\alpha/2}}ds\frac{\sigma^{n-2}}{(1+\sigma)^{n+1}}d\sigma\label{b15}\\
&=&C\int_0^Rr^{1-\alpha/2}\left(\frac{r}{\sqrt{r^2-\sigma^2}}-1\right)^{1-\alpha/2}\frac{\sigma^{n-2}}{(1+\sigma)^{n+1}}d\sigma\nonumber\\
&=&C\int_0^Rr^{1-\alpha/2}\left(\frac{\sigma^2}{\left(r+\sqrt{r^2-\sigma^2}\right)\sqrt{r^2-\sigma^2}}\right)^{1-\alpha/2}\frac{\sigma^{n-2}}{(1+\sigma)^{n+1}}d\sigma\nonumber\\
&\leq&C\int_0^Rr^{1-\alpha/2}\left(\frac{1}{r^2}\right)^{1-\alpha/2}\frac{\sigma^{n-2+2-\alpha}}{(1+\sigma)^{n+1}}d\sigma\nonumber\\
&=&Cr^{\alpha/2-1}\int_0^R\frac{\sigma^{n-\alpha}}{(1+\sigma)^{n+1}}d\sigma\nonumber\\
&\leq&Cr^{\alpha/2-1}\nonumber\\
&=&\frac{C}{r^{1-\alpha/2}}.\label{b16}
\end{eqnarray}

In the above, we derived (\ref{b13}) by letting ~$|y'|=\sigma$~. (\ref{b14}) is valid because
\begin{eqnarray}
\frac{1}{(s^2-1)^{\alpha/2}}&=&\frac{1}{(s+1)^{\alpha/2}}\frac{1}{(s-1)^{\alpha/2}}\nonumber\\&\leq&\frac{1}{(1+1)^{\alpha/2}}\frac{1}{(s-1)^{\alpha/2}}\nonumber\\&=&\frac{1}{2^{\alpha/2}}\frac{1}{(s-1)^{\alpha/2}}\nonumber\\&\leq&\frac{1}{(s-1)^{\alpha/2}}.\nonumber
\end{eqnarray}
Since ~$R$~ is fixed and ~$\sigma^2\leq R^2$~, when~$r$~ is sufficiently large ( much larger than ~$R$~),~we have~$r^2-\sigma^2>0$, and the value of ~$(r^2-\sigma^2)^{\frac{1-\alpha}{2}}$~ can be dominated by ~$(r^2)^{\frac{1-\alpha}{2}}~$(i.e. ~$r^{1-\alpha}$), this verifies (\ref{b15}).

For the ~$\epsilon>0$~ given above and the fixed ~$R$~, since ~$0<\alpha<2$~, then by (\ref{b16}) we can easily get
\begin{equation}
I_{22}\leq C\frac{1}{r^{1-\alpha/2}}<\epsilon,\label{b17}
\end{equation}
for sufficiently large ~$r$~.

From (\ref{b3}), (\ref{b0}), (\ref{b4}), (\ref{b5}), and(\ref{b17}), we derive
\begin{equation}
\left|\frac{\partial u}{\partial x_i}(x)\right|<C\epsilon,
\end{equation}
for sufficiently large ~$R$~ and much larger $r$.

The fact that ~$\epsilon$~ is arbitrary implies
\begin{equation}
\left|\frac{\partial u}{\partial x_i}(x)\right| = 0.
\end{equation}
This proves (\ref{b1}).
\medskip

Now, let's prove (\ref{b2}). Similarly, for fixed ~$x\in B_r(x_r)\subset\mathbb{R}^n_+$~, through an elementary calculation, one can derive that
\begin{eqnarray}
\frac{\partial u}{\partial x_n}(x)&=&\int_{|y-x_r|>r}\left(\frac{\alpha (r- x_n)}{r^2-|x-x_r|^2}+\frac{n(y_n-x_n)}{|y-x|^2}\right)P_r(x-x_r,y-x_r)u(y)dy\nonumber\\
&=&\int_{|y-x_r|>r}\frac{\alpha (r- x_n)}{r^2-|x-x_r|^2}P_r(x-x_r,y-x_r)u(y)dy\nonumber\\&\quad&+\int_{|y-x_r|>r}\frac{n(y_n-x_n)}{|y-x|^2}P_r(x-x_r,y-x_r)u(y)dy\nonumber\\
&:=&J_1+J_2.\label{cc1}
\end{eqnarray}

Similarly to ~$I_2$, for sufficiently large ~$r$, we can also derive
\begin{equation}
|J_2|\leq C\epsilon,
\end{equation}
for any ~$\epsilon>0$~. That is
\begin{equation}
J_2\rightarrow0, \quad\text{as}~~r\rightarrow\infty.\label{cc2}
\end{equation}

Now we estimate ~$J_1$.
\begin{eqnarray*}
J_1&=&\frac{\alpha (r- x_n)}{r^2-|x-x_r|^2}\int_{|y-x_r|>r}P_r(x-x_r,y-x_r)u(y)dy\nonumber\\
&=&\frac{\alpha (r- x_n)}{2x_nr-|x|^2}u(x).
\end{eqnarray*}
It follows that
\begin{equation}
J_1\rightarrow\frac{\alpha}{2x_n}u(x),\quad \text{as}~~r\rightarrow\infty.\label{cc4}
\end{equation}

By (\ref{cc1}), (\ref{cc2}), and(\ref{cc4}), for each fixed ~$x\in B_r(x_r)\subset\mathbb{R}^n_+$, letting ~$r\rightarrow\infty$~, we arrive at
$$
\frac{\partial u}{\partial x_n}(x)=\frac{\alpha}{2x_n}u(x).
$$
This verifies (\ref{b2}), and hence completes the proof of Theorem \ref{mthm1}.\\

\section{$\hat{u}(x)$~is ~$\alpha$~-harmonic in ~$B_r(x_r)$~}
In this section, we prove
 \begin{thm} $\hat{u}(x)$ defined by (\ref{max}) in the previous section is $\alpha$-harmonic in $B_r(x_{r})$.
 \label{thm3.1}
 \end{thm}

The proof consists of two parts. First we show that $\hat{u}$ is harmonic in the average sense (Lemma \ref{lem3.1}), then we show that it is $\alpha$-harmonic (Lemma \ref{lem3.2}).

Let
\begin{equation}
\varepsilon^{(r)}_{\alpha}(x)=\left\{\begin{array}{ll}
0, &|x|<r.\\
\frac{\Gamma(n/2)}{\pi^{\frac{n}{2}+1}} \sin\frac{\pi\alpha}{2}\frac{r^\alpha}{(|x|^2-r^2)^{\frac{\alpha}{2}}|x|^n}, &|x|>r.
\end{array}
\right.
\end{equation}

We say that $u$ is $\alpha$-harmonic in the average sense (see \cite{L}) if for small $r$,
$$\varepsilon^{(r)}_{\alpha}\ast u(x)=u(x).$$

Let
\begin{eqnarray}
&&P_{r}(x-x_r,y-x_r)\nonumber\\&=&\left\{\begin{array}{ll}
\frac{\Gamma(n/2)}{\pi^{\frac{n}{2}+1}} \sin\frac{\pi\alpha}{2}
\left[\frac{r^{2}-|x-x_r|^{2}}{|y-x_r|^{2}-r^{2}}\right]^{\frac{\alpha}{2}}\frac{1}{|x-y|^{n}},\qquad& |y-x_r|>r,|x-x_r|<r\\
0,& \text{elsewhere}.
\end{array}
\right.
\end{eqnarray}

\begin{lem}
 Let $u(x)$ be any measurable function outside $B_r(x_{r})$ for which
 \begin{equation}
 \int_{\mathbb{R}^n}\frac{|u(z)|}{(1+|z-x_{r}|)^{n+\alpha}}dz<\infty.
 \label{d1}
\end{equation}
 Let \begin{equation}
\hat{u}(x)=\left\{\begin{array}{ll}
\int_{|y-x_{r}|>r}P_{r}(y-x_{r},x-x_{r})u(y)dy,& |x-x_{r}|<r,\\
u(x),& |x-x_{r}|\geq r.
\end{array}
\right.
\end{equation}
Then $\hat{u}(x)$ is $\alpha$-harmonic in the average sense in $B_{r}(x_{r})$, i.e. for
sufficiently small $\delta$, we have
\begin{equation}
(\varepsilon_\alpha^{(\delta)}\ast\hat{u})(x)=\hat{u}(x),\quad   |x-x_{r}|<r,
\end{equation}
where $\ast$ is the convolution.
\label{lem3.1}
\end{lem}

\textbf{Proof.}

The outline is as follows.

\textit{i)} \quad Approximate $u$ by a sequence of smooth, compactly supported functions $\{u_k\}$, such that $u_k(x) \ra u(x)$ and
\begin{equation}
\int_{|z-x_{r}|>r}\frac{|u_k(z)-u(z)|}{|z-x_{r}|^n(|z-x_{r}|^2-r^2)^\frac{\alpha}{2}}dz \ra 0.
\label{d2}
\end{equation}

This is possible under our assumption (\ref{d1}).

\textit{ii)} \quad For each $u_k$, find  a signed measure $\nu_k$ such that $supp\,\nu_k\subset B^c_r(x_{r})$ and $$u_k(x)=U_{\alpha}^{\nu_k}(x), \;\; |x-x_{r}|>r.$$
Then $$\hat{u}_k(x)=U_{\alpha}^{\nu_k}(x), \quad |x-x_{r}|<r.$$

\textit{iii)} \quad It is easy to see that $\hat{u}_k(x)$ is $\alpha$-harmonic in the average sense for $|x-x_{r}|<r$.
That is, for each fixed small $\delta>0$,
\begin{equation}
(\varepsilon_\alpha^{(\delta)}\ast\hat{u}_k)(x)=\hat{u}_k(x).
\end{equation}
By showing that as $k \ra \infty$
$$\varepsilon_\alpha^{(\delta)}\ast\hat{u}_k \ra \varepsilon_\alpha^{(\delta)}\ast\hat{u},$$
and
$$\hat{u}_k \ra \hat{u},$$
we arrive at
$$(\varepsilon_\alpha^{(\delta)}\ast\hat{u})(x)=\hat{u}(x),\quad  |x-x_{r}|<r.$$
\bigskip

Now we carry out the details.

\textit{i)} \quad There are several ways to construct such a sequence $\{u_k\}$. One is to use the mollifier.
Let
\begin{equation}
u|_{B_k(x_{r})}(x)=\left\{\begin{array}{ll}
u(x),&|x-x_{r}|<k,\\
0,&|x-x_{r}|\geq k,
\end{array}
\right.
\end{equation}
and
\begin{equation}
J_\epsilon(u|_{B_k(x_{r})})(x)=\int_{\mathbb{R}^n}j_\epsilon(x-y)u|_{B_k(x_{r})}(y)dy.
\end{equation}

For any $\delta>0$, let $k$ be sufficiently large (larger than $r$) such that
\begin{equation}
\int_{|z-x_{r}|\geq k}
\frac{|u(z)|}{|z-x_{r}|^n(|z-x_{r}|^2-r^2)^{\frac{\alpha}{2}}}dz<\frac{\delta}{2}.
\end{equation}
For each such $k$, choose $\epsilon_k$ such that
\begin{equation}
\int_{B_{k+1}\backslash B_r}
\frac{|u_k(z)-u|_{B_k(x_{r})}(z)|}{|z-x_{r}|^n(|z-x_{r}|^2-r^2)^{\frac{\alpha}{2}}}dz
<\frac{\delta}{2},
\end{equation}
where $u_k=J_{\epsilon_k}(u|_{B_k(x_{r})})$. It then follows that
\begin{eqnarray*}
&\quad&\int_{|z-x_{r}|>r}\frac{|u_k(z)-u(z)|}{|z-x_{r}|^n(|z-x_{r}|^2-r^2)^{\frac{\alpha}{2}}}dz\\
&\leq&\int_{B_{k+1}(x_{r})\backslash B_r(x_{r})}\frac{|u_k(z)-u|_{B_k(x_{r})}(z)|+
|u|_{B_k(x_{r})}(z)-u(z)|}{|z-x_{r}|^n(|z-x_{r}|^2-r^2)^{\frac{\alpha}{2}}}dz\\
&&+\int_{|z-x_{r}|>k+1}\frac{|u(z)|}{|z-x_{r}|^n(|z-x_{r}|^2-r^2)^{\frac{\alpha}{2}}}dz\\
&=&\int_{B_{k+1(x_{r})}\backslash B_r(x_{r})}\frac{|u_k(z)-u|_{B_k(x_{r})}(z)|}{|z-x_{r}|^n(|z-x_{r}|^2-r^2)^{\frac{\alpha}{2}}}dz+
\int_{|z-x_{r}|\geq k}\frac{|u(z)|}{|z-x_{r}|^n(|z-x_{r}|^2-r^2)^{\frac{\alpha}{2}}}dz\\
&<&\frac{\delta}{2}+\frac{\delta}{2}=\delta.
\end{eqnarray*}

Therefore, as $k \ra \infty$,
\begin{equation}
\int_{|z-x_{r}|>r}\frac{|u_k(z)-u(z)|}{|z-x_{r}|^n(|z-x_{r}|^2-r^2)^\frac{\alpha}{2}}dz \ra 0.
\label{w1}
\end{equation}

\textit{ii)} \quad For each $u_k$, there exists a signed measure $\psi_k$ such that
\begin{equation}
u_k(x)=U_{\alpha}^{\psi_k}(x).
\end{equation}
 Indeed, let $\psi_k(x)= C (-\lap)^{\frac{\alpha}{2}}u_k(x)$, then
\begin{eqnarray}
U_{\alpha}^{\psi_k}(x)&=&\int_{\mathbb{R}^n}\frac{C}{|x-y|^{n-
 \alpha}}(-\lap)^{\frac{\alpha}{2}}u_k(y)dy\\
 &=&\int_{\mathbb{R}^n}(-\lap)^{\frac{\alpha}{2}}\left[\frac{C}{|x-y|^{n-
 \alpha}}\right]u_k(y)dy\label{d3}\\
 &=&\int_{\mathbb{R}^n}\delta(x-y)u_k(y)dy=u_k(x).
 \end{eqnarray}
Here we have used the fact that $\frac{C}{|x-y|^{n-\alpha}}$ is the fundamental solution of $(-\lap)^{\alpha/2}$.

Let $\psi_k|_{B_r(x_{r})}$ be the restriction of $\psi_k$ on $B_r(x_{r})$ and
\begin{equation}
\tilde{\psi}_k(y)=\int_{|x-x_{r}|<r}P_r(y-x_r, x-x_r)\psi_k|_{B_r(x_{r})}(x)dx,
\end{equation}
we have
 $$U_{\alpha}^{\tilde{\psi}_k}(x)=U_{\alpha}^{\psi_k|_{B_r(x_{r})}}(x), \;\; |x-x_{r}|> r, $$ and
$supp\,\tilde{\psi}_k\subset B^c_r(x_{r}).$ 
Here we use the fact (see (1.6.12$'$) \cite{L}) that
\begin{equation} \frac{1}{|z-x|^{n-\alpha}}=\int_{|y-x_{r}|>r}\frac{P_r(y-x_r,x-x_r)}{|z-y|^{n-\alpha}}dy,\qquad
 \:|x-x_{r}|<r,\: |z-x_{r}|>r.
 \label{Po}
 \end{equation}

 Let $\nu_k=\psi_k-\psi_k|_{B_r(x_{r})}+\tilde{\psi}_k$, then $supp\,\nu_k\subset B^c_r(x_{r})$, and
$$U_{\alpha}^{\nu_k}(x)=U_{\alpha}^{\psi_k}(x)+U_{\alpha}^{\tilde{\psi}_k}(x)
-U_{\alpha}^{\psi_k|_{B_r(x_{r})}}(x)=U_{\alpha}^{\psi_k}(x),\quad |x-x_{r}|>r.$$
That is
$$ u_k(x) = U_{\alpha}^{\nu_k}(x) ,\quad |x-x_{r}|>r.$$

Again by (\ref{Po}), we deduce
$$ \hat{u}_k(x) = U_{\alpha}^{\nu_k}(x) , \quad |x-x_{r}|<r.$$
In this case  $\hat{u}_k$ is $\alpha$-harmonic (in the sense of average) in the region $|x-x_{r}|<r$ (see \cite{L}).

\textit{iii)}\quad For each fixed $x$, we first have
$$\hat{u}_k(x) \ra \hat{u}(x).$$
In fact, by (\ref{w1}),
\begin{eqnarray*}
\hat{u}_k(x)-\hat{u}(x)&=&\int_{|y-x_{r}|>r}P_r(y-x_{r},x-x_{r})[u_k(y)-u(y)]dy\\
&=&C\int_{|y-x_{r}|>r}\frac{(r^2-|x-x_{r}|^2)^{\frac{\alpha}{2}}[u_k(y)-u(y)]}{(
|y-x_{r}|^2-r^2)^{\frac{\alpha}{2}}|x-y|^n}dy\\&\ra&0.
\end{eqnarray*}
 Next, we show that, for each fixed $\delta>0$ and fixed $x$,
 \begin{equation}
 (\varepsilon_\alpha^{(\delta)}\ast\hat{u}_k)(x) \ra (\varepsilon_\alpha^{(\delta)}\ast\hat{u})(x).
 \label{d8}
 \end{equation}

Indeed, 
 \begin{eqnarray*}
 &&(\varepsilon_\alpha^{(\delta)}\ast\hat{u}_k)(x)- (\varepsilon_\alpha^{(\delta)}\ast\hat{u})(x)\\
&=&C\int_{|y-x|>\delta}\frac{\delta^\alpha[\hat{u}_k(y)-
\hat{u}(y)]}{(|x-y|^2-\delta^2)^{\frac{\alpha}{2}}|x-y|^n}dy\\
&=&C\{\int_{\substack{|y-x|>\delta\\ |y-x_{r}|<r-\eta}}
\frac{\delta^\alpha[\hat{u}_k(y)-
\hat{u}(y)]}{(|x-y|^2-\delta^2)^{\frac{\alpha}{2}}|x-y|^n}dy\\
&&+\int_{\substack{|y-x|>\delta \\r-\eta<|y-x_{r}|<r}}
\frac{\delta^\alpha[\hat{u}_k(y)-
\hat{u}(y)]}{(|x-y|^2-\delta^2)^{\frac{\alpha}{2}}|x-y|^n}dy\\
&&+
\int_{\substack{|y-x|>\delta \\|y-x_{r}|>r}}
\frac{\delta^\alpha[\hat{u}_k(y)-
\hat{u}(y)]}{(|x-y|^2-\delta^2)^{\frac{\alpha}{2}}|x-y|^n}dy\}\\
&=&C(I_1+I_2+I_3).
\end{eqnarray*}

For each fixed $x$ with $|x-x_{r}|<r$, choose $\delta$ and $\eta$ such that
$$B_\delta(x)\cap B^c_{r-2\eta}(x_{r})=\emptyset.$$

It follows from (\ref{w1}) that as $k \ra \infty$,
\begin{equation}
I_3=\int_{\substack{|y-x|>\delta \\
                      |y-x_{r}|>r}}\frac{\delta^\alpha[u_k(y)-
u(y)]}{(|x-y|^2-\delta^2)^{\frac{\alpha}{2}}|x-y|^n}dy~ \ra~ 0.
\end{equation}

\begin{eqnarray*}
I_2&=&\int_{\substack{ |y-x|>\delta \\
                      r-\eta<|y-x_{r}|<r}}
\frac{\delta^\alpha \int_{|z-x_{r}|>r}P_r(z-x_{r},y-x_{r})[u_k(z)-u(z)]dz}
{(|x-y|^2-\delta^2)^{\frac{\alpha}{2}}|x-y|^n}dy\\
&=&C\delta^\alpha\int_{|z-x_{r}|>r}\frac{u_k(z)-u(z)}{(|z-x_{r}|^2-r^2)^{\frac{\alpha}{2}}}
\int_{\substack{ |y-x|>\delta \\
                      r-\eta<|y-x_{r}|<r}}
\frac{(r^2-|y-x_{r}|^2)^{\frac{\alpha}{2}}dy}
{(|x-y|^2-\delta^2)^{\frac{\alpha}{2}}|x-y|^n|z-y|^n}dz\\
&=&C\delta^\alpha\int_{|z-x_{r}|>r}\frac{u_k(z)-u(z)}{(|z-x_{r}|^2-r^2)^{\frac{\alpha}{2}}}
\cdot I_{21}(x,z)dz.
\end{eqnarray*}
Noting that in the ring $r-\eta<|y-x_{r}|<r$, we have
$$|x-y|>\eta+\delta.$$
It then follows that
 \begin{eqnarray}
&&I_{21}(x,z)\nonumber\\
&\leq&\frac{1}{(2\eta\delta+\eta^2)^{\frac{\alpha}{2}}(\eta+\delta)^n}
\int_{r-\eta<|y-x_{r}|<r}\frac{(r^2-|y-x_{r}|^2)^{\frac{\alpha}{2}}dy}{|z-y|^n}
\nonumber\\
&=&C\int_{r-\eta}^r (r^2-\tau^2)^{\frac{\alpha}{2}}
\left\{\int_{S_\tau} \frac{1}{|z-y|^n} d\sigma_y\right\}d\tau\nonumber\\
&=&C\int_{r-\eta}^r (r^2-\tau^2)^{\frac{\alpha}{2}}
\left\{\int_0^\pi \frac{\omega_{n-2}(\tau\sin\theta)^{n-2}\tau d\theta}
{(\tau^2+|z-x_{r}|^2-2\tau|z-x_{r}|\cos\theta)^{\frac{n}{2}}}\right\}d\tau\nonumber\\
&=&C\int_{r-\eta}^r (r^2-\tau^2)^{\frac{\alpha}{2}}\frac{1}{\tau^n}
\int_0^\pi \frac{\tau^{n-1}\sin^{n-2}\theta d\theta}{(
(\frac{|z-x_{r}|}{\tau})^2-2\frac{|z-x_{r}|}{\tau}\cos\theta+1)^{\frac{n}{2}}}d\tau
\label{d5}\\
&=&C\int_{r-\eta}^r \frac{(r^2-\tau^2)^{\frac{\alpha}{2}}}{\tau}
\frac{d\tau}{(\frac{|z-x_{r}|}{\tau})^{n-2}((\frac{|z-x_{r}|}{\tau})^2-1)}
\int_{0}^\pi\sin^{n-2}\beta d\beta \label{d6}\\
&<&\frac{Cr^{n-1}}{|z-x_{r}|^{n-2}}\int_{r-\eta}^r
\frac{(r^2-\tau^2)^{\frac{\alpha}{2}}}{|z-x_{r}|^2-\tau^2} d\tau\nonumber\\
&=&\frac{Cr^{n-1}}{|z-x_{r}|^{n-2}}\cdot J.\nonumber
\end{eqnarray}
In the above, to derive (\ref{d6}) from (\ref{d5}), we have made the following substitution (See Appendix in \cite{L}):
$$\frac{\sin\theta}{\sqrt{(\frac{|z-x_{r}|}{\tau})^2-2\frac{|z-x_{r}|}{\tau}\cos\theta+1}}
=\frac{\sin\beta}{\frac{|z-x_{r}|}{\tau}},$$

To estimate the last integral $J$, we consider

\textit{(a)}\quad For $r<|z-x_{r}|<r+1$,
$$J\leq\int_{r-\eta}^r\frac{(r+\tau)^{\frac{\alpha}{2}-1}}{(r-
\tau)^{1-\frac{\alpha}{2}}}d\tau\leq C_{\alpha,r}.$$

\textit{(b)}\quad For $|z-x_{r}|\geq r+1$, obviously,
\begin{equation}
J\sim\frac{1}{|z-x_{r}|^2}, \mbox{ for $|z-x_{r}|$ large}.\nonumber
\end{equation}

In summary,
\begin{equation}
I_{21}(x,z)\sim\left\{\begin{array}{ll}
1,  &\mbox{ for $|z-x_{r}|$ near r},\\
|z-x_{r}|^n,&\mbox{ for $|z-x_{r}|$ large}.
\end{array}
\right.
\nonumber
\end{equation}

Therefore, by (\ref{w1}), as $k \ra \infty$,
\begin{equation}
I_2=\delta^\alpha\int_{|z-x_{r}|>r}\frac{u_k(z)-u(z)}{(|z-x_{r}|^2 -r^2)^{\alpha/2}} I_{21}(x,z)dz ~\ra~ 0.
\end{equation}

Now what remains is to estimate
$$I_1=\delta^\alpha\int_{|z-x_{r}|>r}\frac{u_k(z)-u(z)}{(|z-x_{r}|^2-r^2)^{\frac{\alpha}{2}}}
I_{11}(x,z)dz,$$
where
$$I_{11}(x,z)=\int_{\substack{|y-x|>\delta \\
                      |y-x_{r}|< r-\eta}}
 \frac{(r^2-|y-x_{r}|^2)^{\frac{\alpha}{2}}dy}
 {(|x-y|^2-\delta^2)^{\frac{\alpha}{2}}|x-y|^n|z-y|^n}.$$

 \begin{eqnarray}
 I_{11}(x,z)
&\leq&\frac{r^\alpha}{\delta^n}\int_{\substack{ |y-x|>\delta \\
                      |y-x_{r}|< r-\eta}}
 \frac{dy}{(|x-y|^2-\delta^2)^{\frac{\alpha}{2}}|z-y|^n}\\
&\leq&\frac{r^\alpha}{\delta^n(|z-x_{r}|-r+\eta)^n}
\int_{\delta<|y-x|<2r}\frac{dy}
 {(|x-y|^2-\delta^2)^{\frac{\alpha}{2}}}\\
 &=&\frac{r^\alpha}{\delta^n(|z-x_{r}|-r+\eta)^n}
\int_\delta^{2r}
\frac{\omega_{n-1}\tau^{n-1}d\tau}{(\tau^2-\delta^2)^{\frac{\alpha}{2}}}\\
&\leq&\frac{C}{|z-x_{r}|^n}.
\end{eqnarray}

By (\ref{w1}), as $k \ra \infty$, we have $I_1 \ra 0.$ This verifies (\ref{d8}) and hence completes the proof.
\medskip

\begin{lem}
\begin{equation}
\lim_{r\rightarrow0}\frac{1}{r^\alpha}\left[u(x)-\varepsilon^{(r)}_{\alpha}\ast u(x)\right]
=c(-\lap)^\frac{\alpha}{2}u(x).
\label{c0}
\end{equation}
where $c=\frac{\Gamma(n/2)}{\pi^{\frac{n}{2}+1}} \sin\frac{\pi\alpha}{2}$.
\label{lem3.2}
\end{lem}

\textbf{Proof.}
\begin{eqnarray}
&&\frac{1}{r^\alpha}\left[u(x)-\varepsilon^{(r)}_{\alpha}\ast u(x)\right]\nonumber\\
&=&\frac{1}{r^\alpha}u(x)- c\int_{|y-x|>r}\frac{u(y)}{(|x -y|^2-r^2)^{\frac{\alpha}{2}}|x-y|^n}dy\nonumber\\
&=&c\int_{|y-x|>r}
\frac{u(x)-u(y)}{(|x-y|^2-r^2)^{\frac{\alpha}{2}}|x-y|^n}dy.
\label{c1}
\end{eqnarray}
Here we have used the property that $$\int_{|y-x|>r}\varepsilon^{(r)}_{\alpha}(x-y)=1.$$

 Compare (\ref{c1}) with
$$(-\lap)^\frac{\alpha}{2}u(x)=\lim_{r\rightarrow0}\int_{|y-x|>r}
\frac{u(x)-u(y)}{|x-y|^{\alpha+n}}dy.$$
One may expect that
$$\lim_{r\rightarrow0}\int_{|y-x|>r}
\frac{u(x)-u(y)}{|x-y|^{\alpha+n}}dy=\lim_{r\rightarrow0}\int_{|y-x|>r}
\frac{u(x)-u(y)}{(|x-y|^2-r^2)^{\frac{\alpha}{2}}|x-y|^n}dy.$$

Indeed, consider
\begin{eqnarray}
&&\int_{|y-x|>r}
\frac{u(x)-u(y)}{|x-y|^{n}}\left(\frac{1}{(|x-y|^2-r^2)^{\frac{\alpha}{2}}}
-\frac{1}{|x-y|^{\alpha}}\right)dy\nonumber\\
&=&\int_{r<|y-x|<1}\frac{u(x)-
u(y)}{|x-y|^{n}}\left(\frac{1}{(|x-y|^2-r^2)^{\frac{\alpha}{2}}}
-\frac{1}{|x-y|^{\alpha}}\right)dy\nonumber\\
&&+\int_{|y-x|\geq1}\frac{u(x)-
u(y)}{|x-y|^{n}}\left(\frac{1}{(|x-y|^2-r^2)^{\frac{\alpha}{2}}}
-\frac{1}{|x-y|^{\alpha}}\right)dy\nonumber\\
&=&I_{1}+I_{2}.\label{c2}
\end{eqnarray}

It is easy to see that as $r\rightarrow0$, $I_2$ tends to zero. Actually, same conclusion is true for $I_1$.
\begin{eqnarray}
I_1&=&\int_{r<|y-x|<1}\frac{\nabla u(x)(y-x)+O(|y-x|^2)}{|x-
y|^{n}}\left(\frac{1}{(|x-y|^2-r^2)^{\frac{\alpha}{2}}}
-\frac{1}{|x-y|^{\alpha}}\right)dy\nonumber\\
\label{c2.5}\\
&\leq&C\int_{r<|y-x|<1}\frac{|x-
y|^2}{|x-y|^n}\left(\frac{1}{(|x-y|^2-r^2)^{\frac{\alpha}{2}}}
-\frac{1}{|x-y|^{\alpha}}\right)dy\label{c3}\\
&=&C\int_r^1\frac{\tau^2}{\tau^n}\left(\frac{1}{(\tau^
2-r^2)^{\frac{\alpha}{2}}}
-\frac{1}{\tau^{\alpha}}\right)\tau^{n-1}d\tau\label{c4}\\
&\leq&C\int_1^\infty\left(\frac{1}{r^\alpha(s^
2-1)^{\frac{\alpha}{2}}}
-\frac{1}{r^\alpha s^\alpha}\right)sr^2ds\label{c5}\\
&=&Cr^{2-\alpha}\int_1^\infty\left(\frac{s^\alpha-(s^
2-1)^{\frac{\alpha}{2}}}{(s^
2-1)^{\frac{\alpha}{2}}s^\alpha}\right)sds.\label{c6}
\end{eqnarray}
Equation (\ref{c2.5}) follows from the Taylor expansion. Due to symmetry, we have
$$\int_{r<|y-x|<1}\frac{\nabla u(x)(y-x)}{|x-
y|^{n}}\left(\frac{1}{(|x-y|^2-r^2)^{\frac{\alpha}{2}}}
-\frac{1}{|x-y|^{\alpha}}\right)dy=0$$
 and get (\ref{c3}). By letting $|y-x|=\tau$ and $\tau=rs$ respectively, one obtains (\ref{c4}) and (\ref{c5}).
It is easy to see that the integral in (\ref{c6}) converges near 1. To see that it also converges near infinity, we estimate
$$s^\alpha-(s^2-1)^{\frac{\alpha}{2}}.$$

Let$f(t)=t^{\alpha/2}$. By the \emph{mean value theorem},

\begin{eqnarray*}
f(s^2)-f(s^2-1)
&=&f^\prime(\xi)(s^2-(s^2-1))\\
&=&\frac{\alpha}{2}\xi^{\frac{\alpha}{2}-1}\sim s^{\alpha-2}, \mbox{ for $s$ sufficiently large.}
\end{eqnarray*}

This implies that
$$\frac{s^\alpha-(s^2-1)^{\frac{\alpha}{2}}}{(s^2-1)^{\frac{\alpha}{2}}s^\alpha}s
\sim\frac{s^{\alpha-2}s}{(s^2-1)^{\frac{\alpha}{2}}s^\alpha}
\sim\frac{1}{s^{1+\alpha}}.$$
Now it is obvious that (\ref{c6}) converges near infinity. Thus we have
$$\int_1^\infty\left(\frac{s^\alpha-(s^
2-1)^{\frac{\alpha}{2}}}{(s^
2-1)^{\frac{\alpha}{2}}s^\alpha}\right)sds<\infty.$$
Since $0<\alpha<2$, as $r\rightarrow0$, (\ref{c6}) goes to zero, i.e. $I_1$ converges to zero. Together with (\ref{c1}) and (\ref{c2}), we get (\ref{c0}). This proves the lemma.

\bigskip

{\em Authors' Addresses and E-mails:}
\medskip

Wenxiong Chen 

Department of Mathematical Sciences

Yeshiva University

New York, NY, 10033 USA

wchen@yu.edu
\medskip

Congming Li

Department of Applied Mathematics

University of Colorado,

Boulder CO USA

cli@clorado.edu
\medskip 

Lizhi~Zhang

School~of~Mathematics~and~Information~Science

Henan~Normal~University

azhanglz@163.com
\medskip

Tingzhi~Cheng

Department~of~Mathematics

Shanghai~JiaoTong~University

nowitzki1989@126.com.

\end{document}